\documentclass[preprint,12pt]{elsarticle}
\usepackage{amsmath}
\usepackage{amsfonts}
\usepackage{amssymb}
\usepackage{array,longtable}
\usepackage{amscd}
\usepackage[cp1251]{inputenc}
\usepackage[english]{babel}
\usepackage[matrix,arrow,curve]{xy}

\begin{document}

\begin{frontmatter}
\title{Interlacing properties of the eigenvalues of some matrix classes}
\author{O.Y. Kushel}
\ead{kushel@mail.ru}
\address{Institut f\"{u}r Mathematik, MA 4-5, Technische Universit\"{a}t Berlin, \\ D-10623 Berlin, Germany}

\begin{abstract}
We establish the eigenvalue interlacing property (i.e. the smallest real eigenvalue of a matrix is less than the smallest real eigenvalue of any its principal submatrix) for the class of matrices, introduced by Kotelyansky (all principal and all almost principal minors of these matrices are positive). We show that certain generalizations of Kotelyansky and totally positive matrices also possess this property. We prove some interlacing inequalities for the other eigenvalues of Kotelyansky matrices.
\end{abstract}

\end{frontmatter}

\newtheorem{thm}{Theorem}
\newtheorem{lem}{Lemma}
\newtheorem{prop}{Proposition}
\newtheorem{ob}{Observation}
\newtheorem{cor}{Corollary}
\newproof{pf}{Proof}
\newproof{pot5}{Proof of Theorem 5}
\newproof{pot8}{Proof of Theorem 8}
\newproof{pot10}{Proof of Theorem 10}

\section{Introduction}
The common properties of Hermitian positive semidefinite, $M$- and totally positive matrices allow to unite them into one class, called $\tau$-matrices (see \cite{ENS}, \cite{HERMS1}-\cite{HOLS}).
We refer to the following definitions and notations (see \cite{ENS}, p. 155-156).

Let, as usual, $[n] = \{1, \ \ldots, \ n\}$, $\sigma(A)$ denote the spectrum of a $n \times n$ matrix $\mathbf A$ and let $l(A):= \min\{\sigma(A)\cap{\mathbb R}\}$. If $A$ has no real eigenvalues then $l(A) := \infty$. A matrix $\mathbf A$ is called a {\it $\tau$-matrix} if $l(A) \geq 0$ and $$l(A(\alpha)) \leq l(A(\beta)) < \infty \qquad \mbox{whenever} \qquad \emptyset \neq \beta \subseteq \alpha \subseteq [n]. \eqno(1)$$ (Here $A(\alpha)$, $A(\beta)$ denote principal submatrices of $\mathbf A$, formed by all the rows and columns with indices from the set $\alpha$ (respectively, $\beta$)). A matrix $\mathbf A$ is said to be a {\it strictly $\tau$-matrix ($\tau^<$-matrix) } if, in addition, $l(A) > 0$ and inequalities (1) are all strict.

The property of a matrix $\mathbf A$ defined by the following two conditions
\begin{enumerate}
\item[\rm 1.] each principal submatrix of $\mathbf A$ has a real eigenvalue;
\item[\rm 2.] the minimal real eigenvalues of any two submatrices of $\mathbf A$ satisfy inequalities (1);
\end{enumerate}
is called an {\it eigenvalue monotonicity}.

We study a class of matrices introduced by Kotelyansky (see \cite{KOT1}). These matrices (later called SK-matrices) are defined by the following conditions:
\begin{enumerate}{}
\item[\rm 1.] all of their principal minors are positive;
\item[\rm 2.] all of their minors obtained from principal minors by deleting one row and one
    column with different indices are also positive.
\end{enumerate}

We generalize this definition in the following way.

A $n \times n$ matrix $\mathbf A$ is called {\it strictly J-sign-symmetric Kotelyansky (SJSK)} if it satisfies the following conditions.
\begin{enumerate}
\item[\rm 1.] $\det(A) > 0$.
\item[\rm 2.] If we take any $k \times k$ $(k = 2, \ \ldots, \ n)$ principal submatrix of $\mathbf A$, calculate all its minors of order $k-1$ and arrange them into a matrix in the lexicographic order, then the obtained $k \times k$ matrix will be diagonally similar to a positive matrix.
\end{enumerate}

We show that Kotelyansky matrices as well as mentioned above generalizations satisfy the eigenvalue monotonicity property. Our main results we state in the following two theorems.

\setcounter{thm}{4}
\begin{thm}
An SK-matrix $\mathbf A$ is a $\tau^<$-matrix.
\end{thm}
\setcounter{thm}{7}
\begin{thm}
An SJSK-matrix $\mathbf A$ is a $\tau^<$-matrix.
\end{thm}

In addition, we show that SK-matrices satisfy certain eigenvalue interlacing properties.

\setcounter{thm}{9}
\begin{thm} Let $\mathbf A$ be an $n \times n$ SK matrix. Then $\mathbf A$ has $n$ positive simple eigenvalues:
$$\rho(A) = \lambda_1 > \lambda_2 > \ldots > \lambda_n > 0. $$
Moreover, if $\mu^{(r)}_1 > \mu^{(r)}_2 > \ldots > \mu^{(r)}_{n-1} > 0$ are the eigenvalues of the principal submatrix of $\mathbf A$ obtained by deleting the $r$th row and the $r$th column, then
$$\lambda_j > \mu_j^{(r)}> \lambda_{j+1}, \qquad j = 1, \ \ldots, \ n-1, $$
for $r = 1$ or $r = n$.
\end{thm}
\setcounter{thm}{0}

The organization of this paper is as follows. In Section 2, we recall some definitions used in the paper and study basic properties of Kotelyansky matrices. We mainly concentrate on geometric properties of a linear operator defined by a Kotelyansky matrix. This properties are to be used in the following sections. We introduce certain generalizations of the class of totally positive and Kotelyansky matrices: these are totally J-sign-symmetric and J-sign-symmetric Kotelyansky matrices, respectively. Our main results concerning eigenvalue monotonicity of strictly Kotelyansky and J-sign-symmetric Kotelyansky matrices are proved in Section 3. In Section 4, we study the interlacing properties of strictly Kotelyansky matrices which is analogous that strictly totally positive matrices have.

\section{Kotelyansky matrices: basic definitions and properties}

\subsection{Preliminary notes}
 For $[n]:= \{1, \ \ldots, \ n\}$, as usual, $A\begin{pmatrix}
  i_1 & \ldots & i_k \\
  j_1 & \ldots & j_k
\end{pmatrix}$ with $(i_1, \ \ldots, \ i_k)$,  $(j_1, \ \ldots, \ j_k) \subseteq [n]$, $1 \leq k \leq n$, denotes a minor of a $n \times n$ matrix ${\mathbf A}$ constructed of the rows with numbers $i_1, \ \ldots, \ i_k$ and columns $j_1, \ \ldots, \ j_k$. Let us recall the following definitions concerning some classes of matrix minors (see, for example, \cite{BAO}).

A minor of the form $A\begin{pmatrix}
  i_1 & \ldots & i_k \\
  i_1 & \ldots & i_k
\end{pmatrix}$, where  $(i_1, \ \ldots, \ i_k) \subseteq [n]$, is called {\it principal}.

A minor of the form $A\begin{pmatrix}
  i_1 & \ldots & i_{r-1} & i_{r+1} & \ldots & i_{k} \\
  i_1 & \ldots & i_{s-1} & i_{s+1} & \ldots & i_{k}
\end{pmatrix}$, where $1 \leq r,s \leq k$, $r \neq s$ is called {\it almost principal}.

For a given ${\mathbf A} = \{a_{ij}\}_{i,j = 1}^n$, two sets of indices $(i_1, \ \ldots, i_k)$ and $(j_1, \ \ldots, j_k)$ from $[n]$, and each $l \in [n]\setminus (i_1, \ \ldots, i_k)$, $r \in [n]\setminus (j_1, \ \ldots, j_k)$  let us define $$b_{lr} = A \begin{pmatrix}i_1 & \ldots & i_k & l \\ j_1 & \ldots & j_k & r \end{pmatrix}.$$ In such notation, we always understand that the set of indices is arranged in natural increasing order. Then the {\it Sylvester's Determinant Identity} holds for the minors of the $(n - k)\times(n-k)$ matrix ${\mathbf B} = \{b_{lr}\}$ (see, for example, \cite{PINK}, p. 3):
$$B \begin{pmatrix}l_1 & \ldots & l_p \\ r_1 & \ldots & r_p \end{pmatrix} =
A \begin{pmatrix}i_1 & \ldots & i_k \\ j_1 & \ldots & j_k \end{pmatrix}^{p-1} A \begin{pmatrix}i_1 & \ldots & i_k & l_1 & \ldots & l_p \\ j_1 & \ldots & j_k & r_1 & \ldots & r_p \end{pmatrix}.$$

Let us recall the following basic definitions and facts (see, for example, \cite{CARLMAR}, \cite{CRAH}). Given a matrix $\mathbf A$, written in the following block form:
          $${\mathbf A} = \begin{pmatrix} {\mathbf A}_{11} & {\mathbf A}_{12} \\ {\mathbf A}_{21} & {\mathbf A}_{22} \end{pmatrix},$$
where ${\mathbf A}_{11}$ is a nonsingular leading principal submatrix of ${\mathbf A}$. Then the {\it Schur complement} $A|A_{11}$ of ${\mathbf A}_{11}$ in $\mathbf A$ is defined as follows:
 $$A|A_{11} = A_{22}-A_{21}A_{11}^{-1}A_{12}.$$

Considering the matrix $\widehat{{\mathbf A}}$ obtained from $\mathbf A$ by a simultaneous permutation of rows and columns, which put ${\mathbf A}_{11}$ into the upper left corner of $\widehat{{\mathbf A}}$, we can define the Schur complement of an arbitrary (not necessarily leading) principal submatrix ${\mathbf A}_{11}$ of $\mathbf A$.

If the principal submatrix ${\mathbf A}_{11}$ is spanned by the basic vectors $e_{i_1}, \ \ldots, \ e_{i_k}$, we have $A|A_{11} = \{c_{lr}\}$, where
 $$c_{lr} = A\begin{pmatrix} i_1 & \ldots & i_k & l \\ i_1 & \ldots & i_k & r \end{pmatrix}/A\begin{pmatrix} i_1 & \ldots & i_k \\ i_1 & \ldots & i_k \end{pmatrix}, \eqno(2)$$
$l$, $r$ run over $[n]\setminus (i_1, \ \ldots, i_k)$ and $[n]\setminus (j_1, \ \ldots, j_k)$, respectively.

The following formula connects the inverse of the Schur complement $A|A_{11}$ and the principal submatrix of ${\mathbf A}^{-1}$:
$$(A|A_{11})^{-1} = {A^{-1}}_{22}, $$
where $\cdot_{11}$ denotes the principal submatrix spanned by $e_{i_1}, \ \ldots, \ e_{i_k}$, $\cdot_{22}$ denotes the principal submatrix spanned by the remaining basic vectors.

Let us recall also the Schur formula:
 $$\det A = \det A_{11} \det (A|A_{11}). \eqno(3)$$

\subsection{Exterior powers of the space ${\mathbb R}^n$}
Let $e_1, \ \ldots, \ e_n$ be an arbitrary basis in ${\mathbb R}^n$. Let $x_1, \ \ldots, \ x_j$ $ \ (2 \leq j \leq n)$ be any vectors in ${\mathbb R}^n$ defined by the coordinates: $x_i = (x_i^1, \ \ldots, \ x_i^n)$, $i = 1, \ \ldots, \ j$. Then the vector in ${\mathbb R}^{\binom{n}j}$ ($\binom{n}j = \frac{n!}{j!(n-j)!}$) with coordinates of the form
$$(x_1 \wedge \ldots \wedge x_j)^{\alpha}:= (\begin{vmatrix} x_1^{i_1} & \ldots & x_j^{i_1} \\
\ldots & \ldots & \ldots \\
x_1^{i_j} & \ldots & x_j^{i_j} \\
 \end{vmatrix})^T,$$ where $\alpha$ is the number of the set of indices $(i_1, \ \ldots,  \ i_j) \subseteq [n]$ in the lexicographic order, is
called an {\it exterior product} of $x_1, \ \ldots, \
x_j$.

We consider the $j$th exterior power $\wedge^j {\mathbb R}^n$ of the space ${\mathbb R}^n$ as the space ${\mathbb R}^{\binom{n}j}$. The set of all exterior
products of the form $e_{i_1} \wedge \ldots \wedge e_{i_j}$, where
$1 \leq i_1 < \ldots < i_j \leq n$ forms a canonical basis in $\wedge^j {\mathbb R}^n$ (see \cite{GLALU}).

 Any linear operator $A$ on ${\mathbb R}^n$ defines the linear operator $\wedge^j A$ on $\wedge^j
{\mathbb R}^n$ by
$$ (\wedge^j A)(x_1 \wedge \ldots \wedge x_j) = Ax_1 \wedge \ldots \wedge Ax_j.$$
The operator $\wedge^j A$ is called the {\it $j$th exterior power} of $A$. It is easy to see, that $\wedge^1 A = A$ and $\wedge ^n A$ is one-dimensional and coincide with $\det A$.

If ${\mathbf A} =
\{a_{ij}\}_{i,j = 1}^n$ is the matrix of $A$ in a basis $e_1, \ \ldots, \ e_n$, then
the matrix of $\wedge^j A$ in the basis
$\{e_{i_1} \wedge \ldots \wedge e_{i_j}\}$, where $1 \leq i_1<
\ldots < i_j \leq n$, equals the $j$th compound matrix $
{\mathbf A}^{(j)}$ of the initial matrix ${\mathbf A}$. Here the $j$th compound matrix $
{\mathbf A}^{(j)}$ consists of all the minors of the $j$th order
$A\begin{pmatrix}
  i_1 &  \ldots & i_j \\
  k_1 & \ldots & k_j \end{pmatrix}$, where $1 \leq i_1<
\ldots < i_j \leq n, \ 1 \leq k_1< \ldots < k_j \leq n$, of the
initial $n \times n$ matrix ${\mathbf A}$, listed in the lexicographic order (see, for example,
\cite{PINK}).

\subsection{Basic properties of K-matrices}
We recall that an $n \times n$ real matrix ${\mathbf A}$ is said to be
{\it positive} {\it (non-negative)} if all its entries $a_{ij}$
are positive (respectively, nonnegative). The following classes of matrices were defined in \cite{KOT1} (see also \cite{BAO}).

A matrix ${\mathbf A}$ is called {\it a Kotelyansky matrix (K-matrix)} if all its principal minors are positive and all its almost principal minors are nonnegative. A matrix ${\mathbf A}$ is called a {\it strictly Kotelyansky matrix (SK-matrix)} if all its principal and almost principal minors are positive.

Note that the definition of a K (SK) matrix implies ${\mathbf A}$ and ${\mathbf A}^{(n-1)}$ to be nonnegative (respectively, positive). Note also that if $\mathbf A$ is a K-matrix then $\det A > 0$, so ${\mathbf A}$ is nonsingular.

{\bf Example 1.} Let us show a $4 \times 4$ SK matrix having negative minors (see \cite{BAO}, p. 68).
 $${\mathbf A} = \begin{pmatrix} 3 & 2 & 1 & 0.6 \\ 2 & 3 & 2 & 1 \\ 1 & 2 & 3 & 2 \\ 0.6 & 1 & 2 & 3\end{pmatrix}.$$
In this case
$${\mathbf A}^{(2)} = \begin{pmatrix} 5 & 4 & 1.8 & 1 & 0.2 & -0.2 \\ 4 & 8 & 5.4 & 4 & 2.8 & 0.2 \\ 1.8 & 5.4 & 8.64 & 3 & 5.4 & 1.8 \\ 1 & 4 & 3 & 5 & 4 & 1 \\ 0.2 & 2.8 & 5.4 & 4 & 8 & 4 \\ -0.2 & 0.2 & 1.8 & 1 & 4 & 5 \end{pmatrix};$$
$${\mathbf A}^{(3)} = \begin{pmatrix} 8 & 6.6 & 2.4 & 1 \\ 6.6 & 13.32 & 8.28 & 2.4 \\ 2.4 & 8.28 & 13.32 & 6.6 \\ 1 & 2.4 & 6.6 & 8 \end{pmatrix};$$
$${\mathbf A}^{(3)} = \det A = 12.6.$$

So the matrix ${\mathbf A}$ is SK but not STP. In this case, the two minors $A\begin{pmatrix}
 1 & 2 \\
 3 & 4 \end{pmatrix}$ and $A\begin{pmatrix}
 3 & 4 \\
 1 & 2 \end{pmatrix}$, which are not principal and not almost principal, are negative.

Let us mention some obvious properties of K and SK matrices.

\begin{prop} Let $\mathbf A$ be a K (SK) matrix. Then the following matrices are K (respectively, SK).
\begin{enumerate}
\item[\rm 1.] ${\mathbf A}^T$ (the transpose of $\mathbf A$);
\item[\rm 2.] ${\mathbf D}{\mathbf A}{\mathbf D}^{-1}$, where $\mathbf D$ is a diagonal matrix with positive diagonal entries;
\item[\rm 3.] ${\mathbf D}{\mathbf A}^{-1}{\mathbf D}^{-1}$, where $\mathbf D$ is a diagonal matrix with diagonal entries alternately $+1$ and $-1$;
\item[\rm 4.] ${\mathbf P}{\mathbf A}{\mathbf P}^{-1}$, where $\mathbf P$ is a reverse permutation matrix;
\item[\rm 5.] any principal submatrix of $\mathbf A$;
\item[\rm 6.] the Schur complement of any principal submatrix of $\mathbf A$.
\end{enumerate}
\end{prop}

 \begin{pf} Properties 1-5 immediately follows from the definitions. Let us prove Property 6 for SK matrices (for K-matrices, it is proved analogically). Let us take an arbitrary principal submatrix ${\mathbf A}_{11}$ and consider its Schur complement $A|A_{11}$. Formulas (2) and (3) imply that $A|A_{11}$ is positive and $\det(A|A_{11}) > 0$. Applying Sylvester's determinant identity to principal and almost principal minors of ${\mathbf B} = A|A_{11} * \det A_{11}$, we get:
 $$B \begin{pmatrix}l_1 & \ldots & l_p \\ l_1 & \ldots & l_p \end{pmatrix} =
A \begin{pmatrix}i_1 & \ldots & i_k \\ i_1 & \ldots & i_k \end{pmatrix}^{p-1} A \begin{pmatrix}i_1 & \ldots & i_k & l_1 & \ldots & l_p \\ i_1 & \ldots & i_k & l_1 & \ldots & l_p \end{pmatrix}.$$
$$B \begin{pmatrix}l_1 & \ldots & l_{\nu -1} & l_{\nu + 1}& \ldots & l_p \\ l_1 & \ldots & l_{\mu -1} & l_{\mu + 1}& \ldots & l_p \end{pmatrix} = $$ $$
A \begin{pmatrix}i_1 & \ldots & i_k \\ i_1 & \ldots & i_k \end{pmatrix}^{p-2} A \begin{pmatrix}i_1 & \ldots & i_k & l_1 & \ldots & l_{\nu -1} & l_{\nu + 1}& \ldots & l_p \\ i_1 & \ldots & i_k & l_1 & \ldots & l_{\mu -1} & l_{\mu + 1}& \ldots & l_p\end{pmatrix}.$$
The principal and almost principal minors of $A|A_{11}$ are expressed through the products of principal and almost principal minors of $\mathbf A$, thus they are positive.
\end{pf}

\subsection{Geometry of K-matrices}
Given the set $(i_1, \ \ldots, \ i_j) \subseteq [n]$, we denote $L(e_{i_1}, \ \ldots, \ e_{i_j})$ the subspace of ${\mathbb R}^n$ spanned by the basic vectors $e_{i_1}, \ \ldots, \ e_{i_j}$. If we need only the dimension of the subspace, we use the notation $L(j)$. Later we consider a vector $x \in L(j)$ either as a vector from ${\mathbb R}^n$ whose coordinates with indices from $[n]\setminus (i_1, \ \ldots, \ i_j)$ are equal to zero, or as a vector from ${\mathbb R}^j$ (it depends on the convenience). It is not difficult to come from the one representation to the other by putting/deleting additional zeroes.

As usual, $A|_{L(j)}$ denotes the restriction of the operator $A: {\mathbb R}^n \rightarrow {\mathbb R}^n$ to the subspace $L(j)$, not necessarily invariant with respect to $A$. Here we define the operator $A|_{L(j)}: L(j) \rightarrow L(j)$ as follows: $$A|_{L(j)}(x) = {\rm pr}_{L(j)}(Ax).$$
We consider the matrix of $A|_{L(j)}$ either as a $j \times j$ principal submatrix ${\mathbf A}(i_1, \ \ldots, \ i_j)$ of a matrix $\mathbf A$ or as an $n \times n$ matrix with the corresponding submatrix equal to ${\mathbf A}(i_1, \ \ldots, \ i_j)$ while the remaining entries are zeroes.

\begin{ob}
A linear operator $A: {\mathbb R}^n \rightarrow {\mathbb R}^n$ is SK if and only if $\det A > 0$ and $\wedge^{j-1}A|_{L(j)}$ are positive for all $j = 2, \ \ldots, \ n$ and all $j$-dimensional basic subspaces $L(j)$.
\end{ob}
\begin{pf} $\Rightarrow$ Let $\mathbf A$ be a SK matrix. The inequality $\det A > 0$ immediately follows from the definition. Since $\mathbf A$ is positive, $A|_{L(2)}$ are also positive for any $2$-dimensional basic subspace $L(2)$. Now let us consider a restriction $A|_{L(e_{i_1}, \ \ldots, \ e_{i_j})}$ of the operator $A$ to a $j$-dimensional basic subspace $L(e_{i_1}, \ \ldots, \ e_{i_j})$. The matrix of $A|_{L(e_{i_1}, \ \ldots, \ e_{i_j})}$ has the following form:
$${\mathbf A}(i_1, \ \ldots, \ i_j) = \begin{pmatrix}
 a_{i_1 i_1} & \ldots & a_{i_1 i_j} \\
\ldots & \ldots & \ldots \\
 a_{i_j i_1} & \ldots & a_{i_j i_j}
 \end{pmatrix}.$$
The matrix of $\wedge^{j-1} A|_{L(e_{i_1}, \ \ldots, \ e_{i_j})}$ coincides with the $(j-1)$th compound matrix of the principal submatrix ${\mathbf A}(i_1, \ \ldots, \ i_j)$:
 $${\mathbf A}^{(j-1)}(i_1, \ \ldots, \ i_j) = \begin{pmatrix}
A\begin{pmatrix}
 i_1 & \ldots & i_{j-1} \\
 i_1 & \ldots & i_{j-1} \end{pmatrix} & \ldots & A\begin{pmatrix}
 i_1 & \ldots & i_{j-1} \\
 i_2 & \ldots & i_j  \end{pmatrix} \\
 \ldots & \ldots & \ldots \\
 A\begin{pmatrix}
 i_2 & \ldots & i_{j} \\
 i_1 & \ldots & i_{j-1}  \end{pmatrix} & \ldots & A\begin{pmatrix}
 i_2 & \ldots & i_j \\
 i_2 & \ldots & i_j  \end{pmatrix}
 \end{pmatrix}.$$
Since the matrix ${\mathbf A}^{(j-1)}(i_1, \ \ldots, \ i_j) $ consists of almost principal minors of $\mathbf A$, it is positive.

$\Leftarrow$ Now let all $A|_{L(2)}$ be positive. Then it is obvious that $\mathbf A$ is positive. Let us prove that all principal and almost principal minors of the $(j-1)$th $(j = 3, \ \ldots, \ n)$ order of $\mathbf A$ are positive. Indeed, any minor of the form $A\begin{pmatrix}
  i_1 & \ldots & i_{r-1} & i_{r+1} & \ldots & i_{j} \\
  i_1 & \ldots & i_{s-1} & i_{s+1} & \ldots & i_{j}
\end{pmatrix}$, where $1 \leq r,s \leq j$, is an element of the positive matrix ${\mathbf A}^{(j-1)}(i_1, \ \ldots, \ i_j) $, which is the matrix of $\wedge^{j-1} A|_{L(e_{i_1}, \ \ldots, \ e_{i_j})}$. Thus it is positive. The positivity of the minor of the $n$th order comes from the condition $\det A > 0$. \end{pf}

\begin{ob}
A linear operator $A: {\mathbb R}^n \rightarrow {\mathbb R}^n$ is K if and only if $\wedge^{j-1}A|_{L(j)}$ are nonnegative and $\wedge^{j}A|_{L(j)}$ are positive for all $j = 2, \ \ldots, \ n$ and all $j$-dimensional basic subspaces $L(j)$.
\end{ob}

\begin{pf} $\Rightarrow$ The proof just copies the corresponding reasoning of Observation 1.

$\Leftarrow$ The positivity of $\wedge^{j}A|_{L(j)}$ for all $j = 2, \ \ldots, \ n$ and all $j$-dimensional basic subspaces $L(j)$ is equivalent to the positivity of all principal minors of order $2, \ \ldots, \ n$ of the matrix $\mathbf A$. Let us prove the positivity of all diagonal entries of $\mathbf A$. Since any $2 \times 2$ principal submatrix ${\mathbf A}(i_1, \ i_2) = \begin{pmatrix}a_{i_1 i_1} & a_{i_1 i_2} \\ a_{i_2 i_1} & a_{i_2 i_2}\end{pmatrix}$, $(i_1, \ i_2 = 1, \ \ldots, \ n)$ satisfies the inequalities ${\mathbf A}(i_1, \ i_2) \geq 0$, $\det {\mathbf A}(i_1, \ i_2) = a_{i_1 i_1}a_{i_2 i_2} - a_{i_1 i_2}a_{i_2 i_1} > 0$, we have $a_{i_1 i_1}a_{i_2 i_2} > a_{i_1 i_2}a_{i_2 i_1} \geq 0$. The rest of the proof is analogical to the proof of Observation 1. \end{pf}

\subsection{J-sign-symmetric K-matrices}
 Let $J$ be any subset of $[n] := \{1, 2, \ldots, n\}$. Then $J^c:=[n] \setminus J$ and
$$[n] \times [n] = (J\times J)\cup(J^c \times J^c)\cup (J \times J^c) \cup (J^c \times J)$$ is a partition of $[n] \times [n]$ into four pairwise disjoint subsets.

 A matrix ${\mathbf A} = \{a_{ij}\}_{i,j=1}^n$ is called {\it J-sign-symmetric (JS)} if
$$a_{ij} \geq 0 \quad \mbox{on} \quad (J \times J)\cup (J^c \times J^c);$$
and
$$a_{ij} \leq 0 \quad  \mbox{on} \quad (J \times J^c)\cup (J^c \times J).$$

A matrix ${\mathbf A} = \{a_{ij}\}_{i,j=1}^n$ is called {\it strictly J-sign-symmetric (SJS)} if, in addition, it has no zero entries.

We recall a simple fact that {\it a matrix is J-sign-symmetric (strictly J-sign-symmetric) if and only if it can be represented in the
following form:
$${\mathbf A} = {\mathbf D} \widetilde{{\mathbf A}} {\mathbf
D}^{-1}, $$ where $\widetilde{{\mathbf A}}$ is a
nonnegative (respectively, positive) matrix, ${\mathbf D}$ is a diagonal matrix, which
diagonal entries are equal to $\pm 1$.}

A matrix ${\mathbf A}$ of a linear operator $A:{\mathbb R}^n \rightarrow {\mathbb R}^n$ is called {\it totally J--sign-symmetric (TJS)}, if its $j$-th compound matrices ${\mathbf A}^{(j)}$ are J--sign-symmetric for every $j$ $(j = 1, \ \ldots, \ n)$.

A matrix ${\mathbf A}$ of a linear operator $A:{\mathbb R}^n \rightarrow {\mathbb R}^n$ is called {\it strictly totally J--sign-symmetric (STJS)}, if its $j$-th compound matrices ${\mathbf A}^{(j)}$ are strictly J--sign-symmetric for every $j$ $(j = 1, \ \ldots, \ n)$.

It is easy to see, that the class of totally positive matrices (i.e. matrices whose minors are all nonnegative) belongs to the class of TJS matrices, and the class of STP matrices (i.e. matrices whose minors are all positive) belongs to the class of STJS matrices. Moreover, similarity transformations with permutational and diagonal matrices preserves the classes of TJS and STJS matrices (see \cite{KU}, p. 559, Proposition 30).

However, the class of STJS matrices is wider than all the possible combinations of permutational and diagonal transformations applied to STP matrices.

{\bf Example 2.} Let us consider ${\mathbf A}$
$${\mathbf A} = \begin{pmatrix} 5.6 & 1.2 & 0.7 & 0.5 \\ 6.6 & 6.2 & 4.1 & 8.1 \\ 4.4 & 4.4 & 3.5 & 8\\ 1 & 3.8 & 3.4 & 9 \end{pmatrix}.$$
In this case
$${\mathbf A}^{(2)} = \begin{pmatrix} 26.8 & 18.34 & 42.06 & 0.58 & 6.62 & 3.62 \\
 19.36 & 16.52 & 42.6 & 1.12 & 7.4 & 3.85 \\
 20.08 & 18.34 & 49.9 & 1.42 & 8.9 & 4.6 \\
 1.76 & 5.06 & 17.16 & 3.66 & 13.96 & 4.45 \\
 18.88 & 18.34 & 51.3 & 5.5 & 25.02 & 9.36 \\
 12.32 & 11.46 & 31.6 & 1.66 & 9.2 & 4.3 \\
\end{pmatrix}$$

$${\mathbf A}^{(3)} = \begin{pmatrix}
 15.656 & 58.464 & 15.438 & -2.602 \\
 22.008 & 87.992 & 25.676 & -3.532 \\
 4.168 & 19.76 & 7.69 & -0.45 \\
 -9.584 & -35.408 & -8.354 & 2.386 \\
\end{pmatrix};$$

$${\mathbf A}^{(4)} = \det(A) = 3.3928.$$

The matrix $\mathbf A$ is STJS. It is easy to see that the matrix $\widetilde{{\mathbf A}} = {\mathbf P}{\mathbf A}{\mathbf P}^{-1}$ is still positive for any permutation matrix $\mathbf P$, but the matrix $\widetilde{{\mathbf A}}^{(2)}$ is not positive for any permutation matrix $\mathbf P$, except the identity permutation and the reverse permutation. However, the matrix $\widetilde{{\mathbf A}}^{(3)}$ is not positive for both the identity permutation and the reverse permutation.

Let us provide the following generalizations of the classes of K- and SK-matrices.

 A $n \times n$ matrix $\mathbf A$ is called {\it J-sign-symmetric Kotelyansky (JSK)} if it satisfies the following conditions.
\begin{enumerate}
\item[\rm 1.] All the principal minors of $\mathbf A$ are positive.
\item[\rm 2.] For all $j = 2, \ \ldots, \ n$ and all $j \times j$ principal submatrices ${\mathbf A}(j)$ of $\mathbf A$, the $(j - 1)$th compound matrix ${\mathbf A}^{(j-1)}(j)$ of ${\mathbf A}(j)$ is J-sign-symmetric.
\end{enumerate}

 A $n \times n$ matrix $\mathbf A$ is called {\it strictly J-sign-symmetric Kotelyansky (SJSK)} if it satisfies the following conditions.
\begin{enumerate}
\item[\rm 1.] $\det(A) > 0$.
\item[\rm 2.] For all $j = 2, \ \ldots, \ n$ and all $j \times j$ principal submatrices ${\mathbf A}(j)$ of $\mathbf A$, the $(j - 1)$th compound matrix ${\mathbf A}^{(j-1)}(j)$ of ${\mathbf A}(j)$ is strictly J-sign-symmetric.
\end{enumerate}

{\bf Example 3.} Let us examine a $4 \times 4$ matrix $\mathbf A$.

$${\mathbf A} = \begin{pmatrix} 7 & 8 & 7 & 3 \\ 4.5 & 8 & 4 & 4 \\ 3.5 & 3 & 9 & 4 \\ 2.5 & 5.5 & 8 & 7
\end{pmatrix}.$$
In this case, we have
$${\mathbf A}^{(2)} = \begin{pmatrix} 20 & -3.5 & 14.5 & -24 & 8 & 16 \\ -7 & 38.5 & 17.5 & 51 & 23 & 1 \\ 18.5 & 38.5 & 41.5 & 25.5 & 39.5 & 25 \\ -14.5 & 26.5 & 4 & 60 &
20 & -20 \\ 4.75 & 26 & 21.5 & 42 & 34 & -4 \\ 11.75 & 5.5 &
  14.5 & -25.5 & -1 & 31 \end{pmatrix};$$

$${\mathbf A}^{(3)} = \begin{pmatrix} 106.5 & 64.5 & -88.5 & -120 \\ 119.25 &
  80.25 & -100.5 & -144 \\ -140.25 & -87.75 & 132 &
  178.5 \\ -111.75 & -73.5 & 103.5 & 150 \end{pmatrix};$$

$${\mathbf A}^{(4)} = \det(A) = 42.75.$$

As we see, $\mathbf A$ is positive, ${\mathbf A}^{(3)}$ is SJS and $\det(A) > 0$.

For $3 \times 3$ submatrices of $\mathbf A$, we have:

$${\mathbf A}^{(2)}(1,2,3) = \begin{pmatrix} 20 & -3.5  & -24  \\ -7 & 38.5 & 51 \\ -14.5 & 26.5 & 60 \end{pmatrix}; $$
$${\mathbf A}^{(2)}(1,2,4) = \begin{pmatrix} 20  & 14.5  & 8   \\ 18.5  & 41.5  & 39.5   \\ 4.75 & 21.5 & 34 \end{pmatrix}; $$
$${\mathbf A}^{(2)}(1,3,4) = \begin{pmatrix}   38.5 & 17.5  & 1 \\  38.5 & 41.5  & 25 \\  5.5 &
14.5 & 31 \end{pmatrix}; $$
$${\mathbf A}^{(2)}(2,3,4) = \begin{pmatrix} 60 & 20 & -20 \\  42 & 34 & -4 \\  -25.5 & -1 & 31 \end{pmatrix}. $$

Obviously, all of them are SJS. Thus $\mathbf A$ is SJSK.

Note, that {\it if $\mathbf A$ is a K (SK) matrix, then the matrix $$\widetilde{{\mathbf A}} = {\mathbf P}{\mathbf D}{\mathbf A}{\mathbf D}^{-1}{\mathbf P}^{-1}$$ is JSK (respectively, SJSK) for every permutation matrix ${\mathbf P}$ and every nonsingular diagonal matrix ${\mathbf D}$.}
However, Example 2 shows that not every SJSK matrix can be transformed to an SK matrix using permutation and diagonal matrices.

Later we show that SJSK matrices preserve some common spectral properties of SK and STP matrices.

\section{Weak interlacing properties of K and JSK matrices}

\subsection{Spectral radius}
Let us recall some basic statements concerning positive operators. First, let us state the well-known Perron theorem (see, for example, \cite{GANT1}). Here an {\it eigenfunctional} of the operator $A$ is defined as an eigenvector of the adjoint operator $A^*$.

\begin{thm}[Perron] Let the matrix $\mathbf A$ of a linear operator $A : {\mathbb R}^n
\rightarrow {\mathbb R}^n$ be entry-wise positive. Then:
\begin{enumerate}
\item[\rm 1.] The spectral radius $\rho(A)$ is a simple positive
eigenvalue of the operator $A$ different in absolute value from the
remaining eigenvalues.
\item[\rm 2.] The eigenvector $x_1$ corresponding to
the eigenvalue $\lambda_1 = \rho(A)$ is positive.
\item[\rm 3.] The eigenfunctional $x_1^*$ corresponding to
the eigenvalue $\lambda_1 = \rho(A)$ is positive.
\end{enumerate}
\end{thm}

Let us also recall the following upper and lower estimates for the spectral radius $\rho(A)$ (see \cite{BERPL}, p. 28, Corollary 1.6 and  \cite{BERPL}, p. 15, Theorem 3.31).

\begin{thm}
Let $\mathbf B$ be a principal submatrix of a positive matrix $\mathbf A$. Then $\rho(B) \leq \rho(A)$.
\end{thm}

\begin{thm} Let a matrix $\mathbf A$ be positive. If $0 \neq ({\mathbf A}x - \alpha x) \geq 0$ entry-wise for some nonnegative vector $x$, $\alpha > 0$, then $\rho(A) > \alpha$.
\end{thm}

\subsection{SK matrices are $\tau^<$ matrices}

It has been observed (see, for example, \cite{HERMS1}, \cite{ENS}) that the following matrices belong to the class of $\tau$-matrices.
\begin{enumerate}
\item[\rm 1.] Positive semi-definite matrices (in this case, eigenvalue monotonicity follows from the well-known Cauchy interlacing theorem (see, for example, \cite{BE})).
\item[\rm 2.] M-matrices (eigenvalue monotonicity of M-matrices follows from Frobenius results (see \cite{GANT1}).
\item[\rm 3.] Totally positive matrices (follows of Friedland's results (see \cite{FRIED}, p. 99, Corollary)).
\end{enumerate}

Now we spread the results of \cite{FRIED} to the class of SK-matrices.

For simplicity, denote ${\mathbf A}_k$ an $(n-1) \times (n-1)$ principal submatrix of an $n \times n$ matrix $\mathbf A$, obtained by deleting the $k$th row and the $k$th column of $\mathbf A$. Let $\lambda_1, \ \ldots, \ \lambda_n$ be all the eigenvalues of $\mathbf A$, listed in decreasing order of their absolute values, taking into account their multiplicities. Respectively, let $\mu_1^{(k)}, \ \ldots, \ \mu_{n-1}^{(k)}$ be all the eigenvalues of ${\mathbf A}_k$.

\begin{thm} Let $\mathbf A$ be a SK-matrix. Then the following inequalities hold
$$\lambda_1 > \mu_1^{(k)}, \qquad k = 1, \ \ldots, \ n. \eqno(4) $$
$$\mu_{n-1}^{(k)} > \lambda_n, \qquad k = 1, \ \ldots, \ n. \eqno(5) $$
\end{thm}

\begin{pf} The proof literarily repeats the reasoning of \cite{FRIED} (see \cite{FRIED}, p. 98, the proof of Theorem 2). Let us prove (4). Since ${\mathbf A}_k$ is a principal submatrix of ${\mathbf A}$, we apply Theorem 2 and obtain the inequality $\rho(A_k) \leq \rho(A)$. Assume that $\rho(A_k) = \rho(A)$. Since ${\mathbf A}_k$ is positive (considered as $(n-1) \times (n-1)$ matrix), there is $x_1^k > 0$ such that $ {\mathbf A}_kx_1^k = \rho(A_k)x_1^k$. Examine the vector $\widetilde{x}_1^k$, defined as follows:
$$\widetilde{x}_1^k = ((x_1^k)_1, \ \ldots, \ (x_1^k)_{k-1}, \ 0, \ (x_1^k)_k, \ \ldots, \ (x_1^k)_{n-1}). $$
Thus ${\mathbf A}\widetilde{x}_1^k = \rho(A_k)\widetilde{x}_1^k + k$, where $k = (\underbrace{0, \ \ldots, \ 0}_{k-1}, \ \sum\limits_{i \neq k}a_{ki}(x_1^k)_i, \ 0, \ \ldots, \ 0)$. Since $\mathbf A$ is positive, we have ${\mathbf A}\widetilde{x}_1^k \geq \rho(A_k)\widetilde{x}_1^k$
entrywise with one strict coordinate inequality, for a nonnegative vector $\widetilde{x}_1^k$ and a positive number $\rho(A_k)$. Applying Theorem 3, we get $\rho(A) > \rho(A_k)$.

Now let us prove (5). Applying Observation 1, we obtain that $\det A > 0$ and $\wedge^{j-1}A|_{L(j)}$ are positive for all $j = 2, \ \ldots, \ n$ and all $j$-dimensional basic subspaces $L_j$. Thus $\wedge^{n-2}A|_{L(n-1)} = \wedge^{n-2}A_k > 0$ for any $k = 1, \ \ldots, \ n$. Using well-known relations between $\wedge^{n-1} A$ and $A^{-1}$, we conclude that ${\mathbf B}_k^{-1} > 0$, where ${\mathbf B}_k = {\mathbf D}_{n-1}{\mathbf A}_k{\mathbf D}_{n-1}^{-1}$ and ${\mathbf D}_{n-1} = {\rm diag}(1, \ -1, \ \ldots, \ (-1)^{n-2})$. Applying the Perron theorem (Theorem 1) to  ${\mathbf B}_k^{-1}$, we obtain that $\rho(B_k^{-1})$ is a simple positive eigenvalue of ${\mathbf B}_k^{-1}$ with the corresponding positive eigenvector $y^k$. Observe that $\rho(B_k^{-1})$ is equal to $\dfrac{1}{\mu_{n-1}^k}$, where $\{\mu_1^k, \ \ldots, \ \mu_{n-1}^k\}$ are the eigenvalues of ${\mathbf A}_k$, listed in the decreasing order of their absolute values. Thus $\mu_{n-1}^k > 0$ and simple. The corresponding eigenvector of ${\mathbf A}_k$ is $x^k = {\mathbf D}_{n-1}y^k$. Put $\widetilde{y}_k$ as follows: $$\widetilde{y}^k = (\widetilde{y}_1^k, 0, -\widetilde{y}_2^k), $$
where $\widetilde{y}_1^k = (y_1^k, \ \ldots, \ y_{k-1}^k)$, $\widetilde{y}_2^k = (y_k^k, \ \ldots, \ y_{n-1}^k)$. Examine ${\mathbf B} = {\mathbf D}_n{\mathbf A}{\mathbf D}_n^{-1}$, where ${\mathbf D}_{n} = {\rm diag}(1, \ -1, \ \ldots, \ (-1)^{n-1})$. Since $\wedge^{n-1} A > 0$ and $\det A > 0$ we have ${\mathbf B}^{-1} > 0$. Observe that ${\mathbf D}_n\widetilde{y}^k = \widetilde{x}^k$, where $$\widetilde{x}^k = (x_1^k, \ \ldots, \ x_{k-1}^k, \ 0, \ x_k^k, \ \ldots, \ x_{n-1}^k ).$$ Then we obtain the following equality:
$${\mathbf B}\widetilde{y}_k = {\mathbf D}_n{\mathbf A}{\mathbf D}_n^{-1}\widetilde{y}_k = {\mathbf D}_n{\mathbf A}\widetilde{x}^k = {\mathbf D}_n \begin{pmatrix}\mu_{n-1}^k x_1^k \\ \ldots \\ \mu_{n-1}^k x_{k-1}^k \\ \sum\limits_{i \neq k}a_{ki}x_i^k \\ \mu_{n-1}^k x_k^k \\ \ldots \\ \mu_{n-1}^k x_{n-1}^k \end{pmatrix} = \mu_{n-1}^k\begin{pmatrix} y_1^k \\ \ldots \\ y_{k-1}^k \\ a\\ - y_k^k \\ \ldots \\ - y_{n-1}^k \end{pmatrix}, $$ where $a = \dfrac{(-1)^k}{\mu_{n-1}^k}\sum\limits_{i \neq k}a_{ki}x_i^k$ is some real number. Then $${\mathbf B}^{-1}\begin{pmatrix}\widetilde{y}_1^k \\ a\\ -\widetilde{y}_2^k \end{pmatrix} = \frac{1}{\mu_{n-1}^k}\widetilde{y}_k.$$ Assume that $a<0$. Let us write $${\mathbf B}^{-1} = \begin{pmatrix}{\mathbf B}^{-1}_{11} & {\mathbf B}^{-1}_{12} \\
{\mathbf B}^{-1}_{21} & {\mathbf B}^{-1}_{22} \end{pmatrix},$$ where ${\mathbf B}^{-1}_{11}$ is of order $k-1$ and ${\mathbf B}^{-1}_{22}$ is of order $n-k+1$. Then
$${\mathbf B}^{-1}_{11}\widetilde{y}_1^k =  \dfrac{1}{\mu_{n-1}^k}\widetilde{y}_1^k + {\mathbf B}^{-1}_{12}\begin{pmatrix} -a \\ \widetilde{y}_2^k \end{pmatrix}.$$ Applying Theorem 2 and Theorem 3 we obtain that $0 < \dfrac{1}{\lambda_{n}} = \rho(B^{-1}) \geq \rho(B^{-1}_{11}) > \dfrac{1}{\mu_{n-1}^k}$. Thus $0 <\lambda_n < \mu_{n-1}^k$. The case $a \geq 0$ is treated analogically. Equality (5) is proved. \end{pf}

\begin{pot5} The statement of Theorem 5 follows from Theorem 4 applied to each principal submatrix of $\mathbf A$. (According to Proposition 1 each principal submatrix of an SK-matrix is also an SK-matrix).
\end{pot5}

\subsection{STJS and SJSK matrices are $\tau^<$ matrices}
Considering the exterior powers $\wedge^k {\mathbb R}^n$ and $\wedge^{n-k} {\mathbb R}^n$ $(1 \leq k < n)$ of the finite-dimensional space ${\mathbb R}^n$, we have: $${\rm dim}(\wedge^k {\mathbb R}^n) = {\rm dim}(\wedge^{n-k} {\mathbb R}^n) = \frac{n!}{k!(n-k)!}.$$
All the exterior products of the type $\{e_{i_1}\wedge \ldots \wedge e_{i_{k}}\}$ where $(i_1, \ \ldots, \ i_k) \subset [n]$ of the initial basic vectors $e_1, \
\ldots, \ e_n$ form a basis in $ \wedge^{k} {\mathbb R}^n$, and all the exterior products of the type $\{e_{j_1}\wedge \ldots \wedge e_{j_{n-k}}\}$ where $(j_1, \ \ldots, \ j_{n-k}) \subset [n]$, form a basis in $ \wedge^{n-k} {\mathbb R}^n$. The bijective linear map $J_k: \wedge^k {\mathbb R}^n \rightarrow \wedge^{n-k} {\mathbb R}^n$ is defined on basic vectors as follows:
$$J_k(e_{i_1}\wedge \ldots \wedge e_{i_{k}}) = (-1)^{p+1}e_{j_1}\wedge \ldots \wedge e_{j_{n-k}}, $$
where $(j_1, \ \ldots, \ j_{n-k}) = [n] \setminus (i_1, \ \ldots, \ i_k)$, $p = i_1 + \ldots + i_k$.

Let us examine an invertible linear operator $A: {\mathbb R}^n \rightarrow {\mathbb R}^n$ and its exterior powers $\wedge^k A: \wedge^{k} {\mathbb R}^n \rightarrow \wedge^{k} {\mathbb R}^n$ and $\wedge^{n-k} A: \wedge^{n-k} {\mathbb R}^n \rightarrow \wedge^{n-k} {\mathbb R}^n$. Applying the well-known Laplace expansion theorem (see, for example, \cite{PINK}, p. 4), we get the identities
$$\det(A)I = (\wedge^k A) (J_k)^{-1} (\wedge^{n-k}A) J_k$$
or
$$ (\wedge^k A)^{-1} = \dfrac{1}{\det A}(J_k)^{-1} (\wedge^{n-k}A) J_k \eqno (6)$$
If the operator $A$ is defined by an $n \times n$ matrix $\mathbf A$, we can re-write Equality (6) in terms of compound matrices:
$$({\mathbf A}^{(k)})^{-1} = \dfrac{1}{\det A}{\mathbf D}_{\binom{n}k}{\mathbf P}_k{\mathbf A}^{(n-k)}{\mathbf P}_k{\mathbf D}_{\binom{n}k}, $$
where ${\mathbf P}_k$ is a permutational matrix corresponding to the reverse permutation of $[\binom{n}k]$, ${\mathbf D}_{\binom{n}k}$ is a diagonal matrix with diagonal entries $d_{\alpha\alpha} = (-1)^{p+1}$, $p = i_1 + \ldots + i_k$, $(i_1, \ \ldots, \ i_k)$ is a set of indices from $[n]$ which has number $\alpha$ in the lexicographic ordering.
Later we focus on the cases $k=1$ and $k = 2$. For $k=1$, we obtain the well-known relation between the inverse matrix ${\mathbf A}^{-1}$ and ${\mathbf A}^{(n-1)}$.
$${\mathbf A}^{-1} = \dfrac{1}{\det A}{\mathbf D}_{n}{\mathbf P}_1{\mathbf A}^{(n-1)}{\mathbf P}_1 {\mathbf D}_{n}, \eqno (7) $$
where ${\mathbf D}_{n} = {\rm diag}\{1, \ -1, \ 1, \ \ldots\}$.
For $k = 2$, we have
$$({\mathbf A}^{(2)})^{-1} = \dfrac{1}{\det A}{\mathbf D}_{\binom{n}2}{\mathbf P}_2{\mathbf A}^{(n-2)}{\mathbf P}_2 {\mathbf D}_{\binom{n}2}, $$
or, applying the identity for the minors of the inverse matrix (see \cite{PINK}, p. 4)
$$({\mathbf A}^{-1})^{(2)} = \dfrac{1}{\det A}{\mathbf D}_{\binom{n}2}{\mathbf P}_2{\mathbf A}^{(n-2)}{\mathbf P}_2 {\mathbf D}_{\binom{n}2}. \eqno(8) $$

Now we recall the following lemma.
\begin{lem} Let a matrix $\mathbf A$ be SJS and its second compound matrix ${\mathbf A}^{(2)}$ be also SJS. Then there exists a permutation matrix $\mathbf P$ and a diagonal matrix $\mathbf D$ with diagonal entries equal to $\pm 1$ such that $${\mathbf A} = {\mathbf D}{\mathbf P}\widetilde{{\mathbf A}}{\mathbf P}^T{\mathbf D}^{-1},$$
where the matrix $\widetilde{{\mathbf A}}$ is positive and its second compound matrix $\widetilde{{\mathbf A}}^{(2)}$ is also positive.
\end{lem}
\begin{pf} For the proof, see \cite{KU2}, Corollary 3 and Theorem 10. \end{pf}

Now we can prove the eigenvalue monotonicity of STJS matrices.

\setcounter{thm}{5}
 \begin{thm} Let $\mathbf A$ be an STJS-matrix. Then Inequalities (4) and (5) hold:
$$\lambda_1 > \mu_1^{(k)}, \qquad k = 1, \ \ldots, \ n. \eqno(4) $$
$$\mu_{n-1}^{(k)} > \lambda_n, \qquad k = 1, \ \ldots, \ n. \eqno(5) $$
\end{thm}

 \begin{pf} First prove Inequality (4). Since $\mathbf A$ is SJS, it can be represented in the form $${\mathbf A} = {\mathbf D}\widetilde{{\mathbf A}}{\mathbf D}^{-1},$$ where $\widetilde{{\mathbf A}}$ is a positive matrix, ${\mathbf D}$ is a diagonal matrix. The following equality holds for any $(n-1)\times(n-1)$ principal submatrix ${\mathbf A}_k$ of $\mathbf A$:
$${\mathbf A}_k = {\mathbf D}_k\widetilde{{\mathbf A}}_k{\mathbf D}_k^{-1},$$ where $\widetilde{{\mathbf A}}_k$ and ${\mathbf D}_k$ are the corresponding $(n-1) \times (n-1)$ principal submatrices of $\widetilde{{\mathbf A}}$ and $\mathbf D$, respectively. Since the above transformations preserve spectra of matrices, we have $\lambda_1(A) = \lambda_1(\widetilde{A})$ and $\mu_1(A_k) = \mu_1(\widetilde{A}_k)$. Repeating the first part of the proof of Theorem 4, we obtain $\lambda_1(\widetilde{A}) > \mu_1(\widetilde{A}_k)$. This implies $\lambda_1(A) > \mu_1(A_k)$.

Let us prove Inequality (5). Since $\mathbf A$ is SJS, we have ${\mathbf A}^{(n-2)}$ is SJS, ${\mathbf A}^{(n-1)}$ is SJS and $\det(A) > 0$. Applying Identities (7) and (8), we obtain $${\mathbf A}^{(n-1)} = (\det A){\mathbf D}_{n}{\mathbf P}_1{\mathbf A}^{-1}{\mathbf P}_1 {\mathbf D}_{n} $$
and
$${\mathbf A}^{(n-2)} = (\det A){\mathbf D}_{\binom{n}2}{\mathbf P}_2({\mathbf A}^{-1})^{(2)}{\mathbf P}_2 {\mathbf D}_{\binom{n}2}. $$
Observing that ${\mathbf D}_{n}^{(2)} = - {\mathbf D}_{\binom{n}2}$, ${\mathbf P}_1^{(2)} = - {\mathbf P}_2$ and applying the Cauchy--Binet formula, we obtain
$${\mathbf A}^{(n-2)} = (\det A)({\mathbf D}_n{\mathbf P}_1{\mathbf A}^{-1}{\mathbf P}_1 {\mathbf D}_n)^{(2)} = \dfrac{1}{\det A}({\mathbf A}^{(n-1)})^{(2)}. $$
Since ${\mathbf A}^{(n-1)}$ and ${\mathbf A}^{(n-2)}$ are SJS, we have, that the matrix ${\mathbf B} = {\mathbf D}_{n}{\mathbf P}_1{\mathbf A}^{-1}{\mathbf P}_1{\mathbf D}_{n} $ is SJS as well as its second compound matrix ${\mathbf B}^{(2)}$. Applying Lemma 1 to $\mathbf B$, we obtain that $${\mathbf B} = {\mathbf D}{\mathbf P}\widetilde{{\mathbf B}}{\mathbf P}^T{\mathbf D}^{-1},$$ where $\widetilde{{\mathbf B}}$ is a positive matrix, its second compound matrix $\widetilde{{\mathbf B}}^{(2)}$ is also positive, $\mathbf D$ is a diagonal matrix with diagonal entries equal to $\pm 1$, $\mathbf P$ is a permutation matrix.

Examine the matrix $\widetilde{{\mathbf A}} = {\mathbf S}{\mathbf A}{\mathbf S}^{-1},$ where ${\mathbf S} = {\mathbf P_1}{\mathbf D}_n{\mathbf D}{\mathbf P}{\mathbf D}_n{\mathbf P_1}$. Applying Identities (7) and (8), we get $$\widetilde{{\mathbf A}}^{(n-1)} = (\det A){\mathbf D}_{n}{\mathbf P_1}\widetilde{{\mathbf A}}^{-1}{\mathbf P_1}{\mathbf D}_{n}  = (\det A){\mathbf D}_{n}{\mathbf P_1}{\mathbf S}{\mathbf A}^{-1}{\mathbf S}^{-1}{\mathbf P_1} {\mathbf D}_{n} = $$ $$
= (\det A){\mathbf D}{\mathbf P}{\mathbf D}_n{\mathbf P}_1{\mathbf A}^{-1}{\mathbf P}_1{\mathbf D}_n{\mathbf P}^T{\mathbf D}^{-1} =(\det A) \widetilde{{\mathbf B}}.$$
Analogically,
$$\widetilde{{\mathbf A}}^{(n-2)} = (\det A)({\mathbf D}_{n}{\mathbf P_1}\widetilde{{\mathbf A}}^{-1}{\mathbf P_1} {\mathbf D}_{n})^{(2)} = (\det A) \widetilde{{\mathbf B}}^{(2)}.$$
Since $\widetilde{{\mathbf B}}$ is a positive matrix and $\widetilde{{\mathbf B}}^{(2)}$ is also positive, we have $\widetilde{{\mathbf A}}^{(n-1)}$ is positive and $\widetilde{{\mathbf A}}^{(n-2)}$ is also positive. Repeating the proof of the second part of Theorem 4, we obtain that $\mu_{n-1}(\widetilde{A}_k) > \lambda_n(\widetilde{A})$ for any $k = 1, \ \ldots, \ n$. The matrix $\widetilde{{\mathbf A}}$ is similar to $\mathbf A$ with the composition of a permutation matrix and diagonal matrices with diagonal entries equal to $\pm 1$. It is easy to see, that such a similarity transformation preserves the spectra of each submatrix of the matrix $\mathbf A$. Thus Inequality (5) holds for $\mathbf A$ as well. \end{pf}

\begin{cor}
An STJS-matrix $\mathbf A$ is a $\tau^<$-matrix.
\end{cor}

Now let us consider the class of SJSK matrices.

\begin{thm} Let $\mathbf A$ be an SJSK-matrix. Then Inequalities (4) and (5) hold:
$$\lambda_1 > \mu_1^{(k)}, \qquad k = 1, \ \ldots, \ n. \eqno(4) $$
$$\mu_{n-1}^{(k)} > \lambda_n, \qquad k = 1, \ \ldots, \ n. \eqno(5) $$
\end{thm}

\begin{pf} The proof of Inequality (4) copies the reasoning of Theorem 6. Let us prove Inequality (5). Let us fix an arbitrary $k$ and examine an $(n-1) \times (n-1)$ submatrix ${\mathbf A}_k$. According to the definition of SJSK matrices, ${\mathbf A}_k^{(n-2)}$ is SJS and ${\mathbf A}^{(n-1)}$ is also SJS.

Applying Identities (7) and (8), we obtain $${\mathbf A}^{(n-1)} = (\det A){\mathbf D}_{n}{\mathbf P}_1{\mathbf A}^{-1}{\mathbf P}_1 {\mathbf D}_{n}; $$
$${\mathbf A}^{(n-2)} = (\det A){\mathbf D}_{\binom{n}2}{\mathbf P}_2({\mathbf A}^{-1})^{(2)} {\mathbf P}_2{\mathbf D}_{\binom{n}2}$$
and
$${\mathbf A}^{(n-2)} = (\det A)({\mathbf D}_n {\mathbf P}_1{\mathbf A}^{-1}{\mathbf P}_1 {\mathbf D}_n)^{(2)}. $$

For simplicity, let us assume the matrix ${\mathbf B} = {\mathbf D}_{n}{\mathbf P}_1{\mathbf A}^{-1} {\mathbf P}_1{\mathbf D}_{n}$ to be positive, otherwise we consider the matrix $\widetilde{{\mathbf B}} = {\mathbf D}{\mathbf B}{\mathbf D}^{-1}$, where ${\mathbf D}$ is the required diagonal matrix. Examine the $(n-1)$-dimensional basic subspace $L_k \subset {\mathbb R}^n$, spanned by the set of basic vectors $\{e_1, \ \ldots, \ e_{k-1}, \ e_{k+1}, \ \ldots, \ e_n\}$, and the corresponding exterior basic subspace $\wedge^{n-2}L_k$ of $\wedge^{n-2}{\mathbb R}^n$, spanned by all the possible exterior products $e_{i_1} \wedge \ldots \wedge e_{i_{n-1}}$ of $\{e_1, \ \ldots, \ e_{k-1}, \ e_{k+1}, \ \ldots, \ e_n\}$. Applying the map $J_2^{-1}$, we obtain the corresponding $(n-1)$-dimensional subspace of $\wedge^2{\mathbb R}^n$ which is spanned by exterior products $e_i \wedge e_k$, $i = 1, \ \ldots, \ k-1$ and $e_k \wedge e_i$, $i = k+1, \ \ldots, \ n$. The corresponding submatrix ${\mathbf B}^{(2)}_{\chi}$ of ${\mathbf B}^{(2)} = ({\mathbf D}_{n}{\mathbf P}_1{\mathbf A}^{-1}{\mathbf P}_1 {\mathbf D}_{n})^{(2)}$ is SJS since ${\mathbf A}_k^{(n-2)}$ is SJS. Thus there exists an $(n-1)\times(n-1)$ diagonal matrix $\widehat{{\mathbf D}}$ such that $\widehat{{\mathbf D}}{\mathbf B}^{(2)}_{\chi}\widehat{{\mathbf D}}$ is positive. Regarding $\widehat{{\mathbf D}}$ on the exterior basic vectors $\{e'_i \wedge e'_k\}$, we arrange the indices $\{1, \ \ldots, \ k-1, \ k+1, \ \ldots, \ n\}$ into two groups: $$ \left\{\begin{array}{cc} i \in I_1 & \mbox{if} \ i < k \ \mbox{and} \ \widehat{d}_{ii} = 1 \ \mbox{or} \ i > k \ \mbox{and} \ \widehat{d}_{i-1i-1} = -1;
\\[10pt] i \in I_2 & \mbox{otherwise.}\end{array}\right.$$
Then we define a permutation $\tau$ of the set $[n]$ as follows:
$$\tau = \{I_1'\}\cup \{k\} \cup \{I'_2\},$$
where $I'_1$ and $I'_2$ are the sets $I_1$ and $I_2$ respectively, ordered lexicographically. Examine the matrix $\widetilde{{\mathbf B}} = {\mathbf P}{\mathbf B}{\mathbf P}^T$, where $\mathbf P$ is a permutation matrix, defined by the permutation $\tau$. We have that $\widetilde{{\mathbf B}}$ is positive since $\mathbf B$ is positive and the $(n-1)\times (n-1)$ submatrix $\widetilde{{\mathbf B}}^{(2)}_{\chi}$ of the matrix $\widetilde{{\mathbf B}}^{(2)}$ is also positive. Then examine the matrix ${\widetilde {\mathbf A}} = {\mathbf S}{\mathbf A}{\mathbf S}^{-1}$, where ${\mathbf S} = {\mathbf P}_1{\mathbf D}_n{\mathbf P}{\mathbf D}_n {\mathbf P}_1$. In this case, we have the equalities
$$\widetilde{{\mathbf A}}^{(n-1)} = (\det A) \widetilde{{\mathbf B}},$$
$$\widetilde{{\mathbf A}}^{(n-2)} = (\det A) \widetilde{{\mathbf B}}^{(2)} $$
and the corresponding equality for the compound of the submatrix ${\mathbf A}_k$:
$$\widetilde{{\mathbf A}}_k^{(n-2)} = (\det A)\widetilde{{\mathbf B}}^{(2)}_{\chi}. $$
Thus both $\widetilde{{\mathbf A}}^{(n-1)}$ and $\widetilde{{\mathbf A}}_k^{(n-2)}$ are positive. Repeating the reasoning from the proof of the second part of Theorem 4, we complete the proof.
\end{pf}

\begin{pot8} We observe that each principal submatrix of an SJSK-matrix is also an SJSK-matrix. Thus the statement of Theorem 8 follows from Theorem 7 applied to each principal submatrix of $\mathbf A$.
\end{pot8}

\section{Interlacing properties of the eigenvalues of SK matrices}

Following Kotelyansky, let us use the notation $A_\lambda = A - \lambda I$. Let us determine a sequence $A_\lambda^1, \ \ldots, \ A_\lambda^n$ of enclosed minors of $A_\lambda$. Here $A_\lambda^{k} = A_\lambda \begin{pmatrix} k & k+1 & \ldots & n \\ k & k+1 & \ldots & n \end{pmatrix}$, $k = 1, \ \ldots, \ n$.
Let us define $\alpha, \beta \in {\mathbb R}$ such that $A^k_\alpha > 0$ and $(-1)^k A^k_\beta > 0$ for all $k = 1, \ \ldots, \ n$. Note that the leading coefficient of $A_\lambda^{k}$ is equal to $(-1)^k$. As it follows, $\alpha = - \infty$ and $\beta = + \infty$ satisfy the given conditions. So we can always find sufficiently big finite numbers $\alpha$ and $\beta$.

The following theorem is proved by Kotelyansky (see \cite{KOT1}, p. 164, Theorem 1).

\setcounter{thm}{8}
\begin{thm}[Kotelyansky]
If every minor in the sequence
$$A_\lambda \begin{pmatrix} 1 & 3 & \ldots & n \\ 2 & 3 & \ldots & n \end{pmatrix}, \ A_\lambda \begin{pmatrix} 2 & 4 & \ldots & n \\ 3 & 4 & \ldots & n \end{pmatrix}, \ \ldots, \ A_\lambda \begin{pmatrix} n-1 \\ n \end{pmatrix}$$
has the same sign that the corresponding minor in the sequence
$$A_\lambda \begin{pmatrix} 2 & 3 & \ldots & n \\ 1 & 3 & \ldots & n \end{pmatrix}, \ A_\lambda \begin{pmatrix} 3 & 4 & \ldots & n \\ 2 & 4 & \ldots & n \end{pmatrix}, \ \ldots, \ A_\lambda \begin{pmatrix} n \\ n-1 \end{pmatrix}$$
for all $\lambda \in [\alpha,\beta]$ then all the minors
$$A_\lambda \begin{pmatrix} 1 & 2 & \ldots & n \\ 1 & 2 & \ldots & n \end{pmatrix}, \ A_\lambda \begin{pmatrix} 2 & 3 & \ldots & n \\ 2 & 3 & \ldots & n \end{pmatrix}, \ \ldots, \ A_\lambda \begin{pmatrix} n \\ n \end{pmatrix}$$
have real simple roots lying in $(\alpha,\beta)$. Moreover, the roots of any two consecutive minors  strictly interlace.
\end{thm}

Let us show, that SK matrices has an interlacing property, which is analogous to the interlacing property of STP matrices.

\begin{pot10} Let us consider the case $r = 1$ (the case $r = n$ is considered analogically). The proof literally repeats the reasoning of \cite{KOT1}, p. 164, Proof of Theorem 1 and \cite{PINK}, p. 136, Proof of Proposition 5.4. First let us check for an $n\times n$ SK matrix $\mathbf A$ that $$A_\lambda \begin{pmatrix}1 & \ldots & n-1 \\ 2 & \ldots & n \end{pmatrix} A_\lambda \begin{pmatrix}2 & \ldots & n \\ 1 & \ldots & n-1 \end{pmatrix} > 0$$ for all $\lambda > 0$. Indeed,
let us expand the first minor as polynomial in $\lambda$ (the second minor is expanded analogically).
$$A_\lambda \begin{pmatrix}1 & \ldots & n-1 \\ 2 & \ldots & n \end{pmatrix} = $$ $$ = (-\lambda)^{n-2}A\begin{pmatrix}1 \\ n \end{pmatrix} + (- \lambda)^{n-3}\sum_{i_1 = 2}^{n-1} A\begin{pmatrix}1 & i_1 \\ n & i_1\end{pmatrix} + (- \lambda)^{n-4}\sum_{i_1, i_2 = 2}^{n-1} A\begin{pmatrix}1 & i_1 & i_2 \\ n & i_1 & i_2 \end{pmatrix} + $$ $$+ \ldots + \det A =$$
$$ = (-1)^{n-2}\lambda^{n-2}A\begin{pmatrix}1 \\ n \end{pmatrix} + (-1)^{n-2} \lambda^{n-3}\sum_{i_1 = 2}^{n-1} A\begin{pmatrix}1 & i_1 \\ i_1 & n \end{pmatrix} + $$ $$ + (-1)^{n-2} \lambda^{n-4}\sum_{i_1, i_2 = 2}^{n-1} A\begin{pmatrix}1 & i_1 & i_2 \\ i_1 & i_2 & n \end{pmatrix} + \ldots + \det A.$$
Since all almost principal minors of $\mathbf A$ are positive, all the coefficients of this polynomial has the same sign $(-1)^{n-2}$. So using Descartes' Rule of Signs we obtain that it has no positive roots. Thus it preserves the sign on $(0, + \infty)$.

Since the matrix $\mathbf A$ is SK we have that any its principal submatrix is also SK. Applying the above reasoning consequently to all enclosed principal submatrices obtained from $\mathbf A$ by deleting the first (or the last) row an column we obtain
$$A_\lambda \begin{pmatrix}k & \ldots & n-k+1 \\ k+1 & \ldots & n-k+2 \end{pmatrix}A_\lambda \begin{pmatrix}k+1 & \ldots & n-k+2 \\ k & \ldots & n-k+1 \end{pmatrix} > 0 $$
$$A_\lambda \begin{pmatrix}k & \ldots & n-k \\ k+1 & \ldots & n-k+1 \end{pmatrix}A_\lambda \begin{pmatrix}k+1 & \ldots & n-k+1 \\ k & \ldots & n-k \end{pmatrix} > 0 $$
for all $k = 2, \ \ldots, \ \lceil\frac{n}{2}\rceil $ and all $\lambda \in (0, + \infty)$.

Let us examine the matrix $\widetilde{\mathbf {A}} = {\mathbf P}{\mathbf A}{\mathbf P}^T$, where the permutation matrix $\mathbf P$ is defined by the permutation $\sigma = (1, \ n, \ 2, \ n-1, \ \ldots, \ l)$, $l = \lfloor\frac{n}{2}\rfloor + 1$.

Since $$\widetilde{A}_\lambda \begin{pmatrix}2k-1 & 2k  & \ldots & n \\ 2k-2 & 2k &  \ldots & n \end{pmatrix}\widetilde{A}_\lambda \begin{pmatrix}2k-2 & 2k & \ldots & n \\ 2k-1 & 2k &  \ldots & n \end{pmatrix} = $$ $$ = A_\lambda \begin{pmatrix}k & \ldots & n-k+1 \\ k+1 & \ldots & n-k+2 \end{pmatrix}A_\lambda \begin{pmatrix}k+1 & \ldots & n-k+2 \\ k & \ldots & n-k+1 \end{pmatrix} > 0;$$
$$ \widetilde{A}_\lambda \begin{pmatrix}2k & 2k+1  & \ldots & n \\ 2k-1 & 2k+1 &  \ldots & n \end{pmatrix}\widetilde{A}_\lambda \begin{pmatrix}2k-1 & 2k+1 & \ldots & n \\ 2k & 2k+1 &  \ldots & n \end{pmatrix} =$$ $$= A_\lambda \begin{pmatrix}k & \ldots & n-k \\ k+1 & \ldots & n-k+1 \end{pmatrix}A_\lambda \begin{pmatrix}k+1 & \ldots & n-k+1 \\ k & \ldots & n-k \end{pmatrix} > 0, $$
for all $k = 1, \ \ldots, \ \lceil\frac{n}{2}\rceil $ and all $\lambda \in (0, + \infty)$,
we obtain that the matrix $\widetilde{\mathbf {A}}$ satisfies the conditions of Theorem 9.
Applying Theorem 9, we complete the proof.
\end{pot10}

This property does not holds for all $r$ even if all the minors of $\mathbf A$ are positive. The counterexample can be found in \cite{PINK}, p. 141.

\end{document}